# The Four Hagge Circles

**Christopher J Bradley**

**Abstract:** The four Hagge circles of the triangles BCD, ACD, ABD, ABC of a cyclic quadrilateral ABCD with respect to an appropriate common axis of inverse similarity share many features which are analysed.

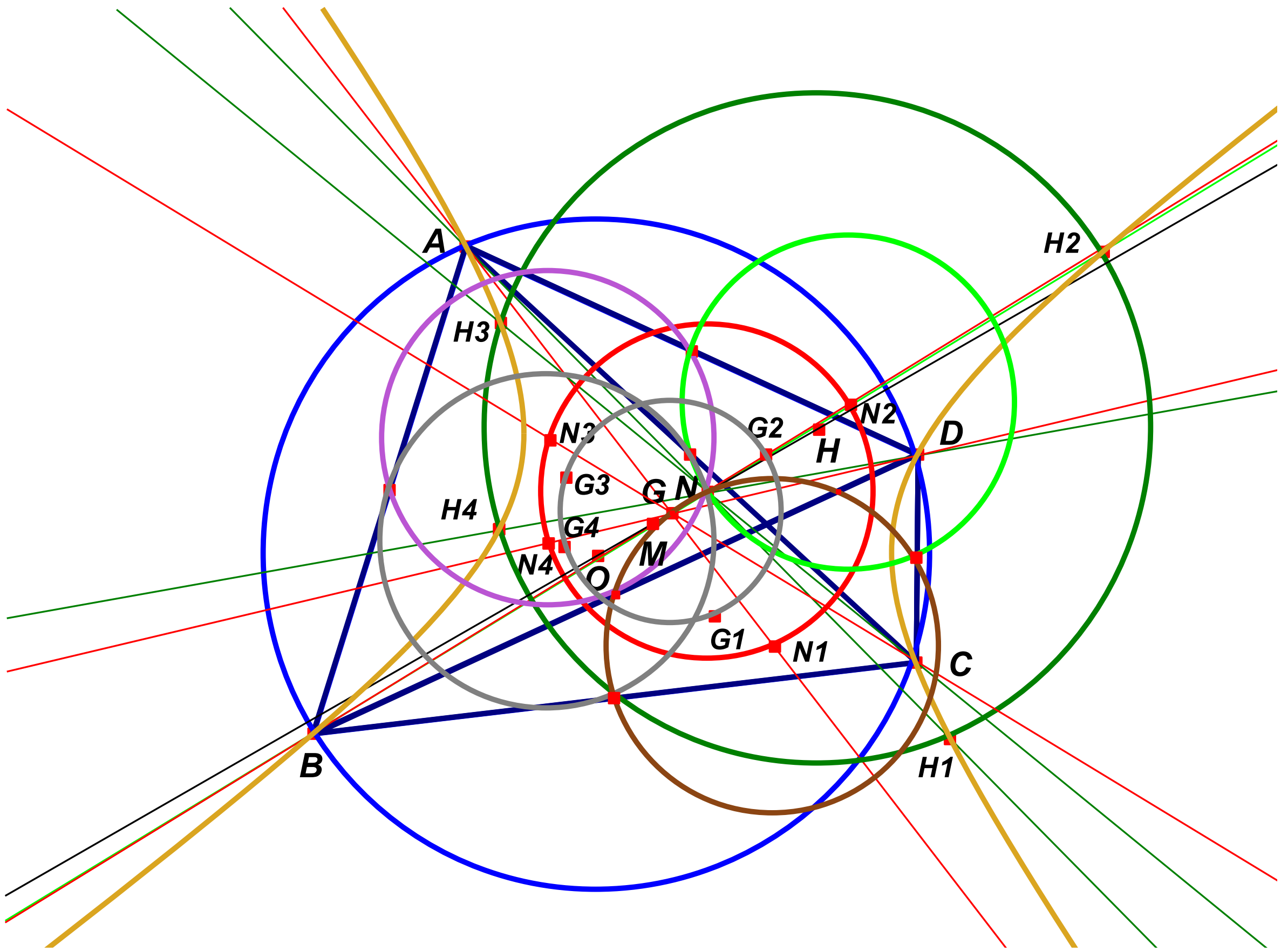

**Fig.1**
**The key points involved in describing the four triangles of a cyclic quadrilateral**

## 1. Preliminary analysis

In later sections we state and prove a number of theorems about the four Hagge circles, with respect to an appropriately chosen centre of inverse similarity *P*, of the four triangles *BCD*, *ACD*, *ABD*, *ABC* that comprise a cyclic quadrilateral *ABCD*. It is desirable therefore to give a description of the basic configuration involved and to prove some straightforward facts that are required later.

We use vectors with origin *O*, the centre of the circumscribing circle, and with *A*, *B*, *C*, *D* having vector positions **a**, **b**, **c**, **d** each of which has magnitude *R*, the radius of the circle. The triangles *BCD*, *ACD*, *ABD*, *ABC* are denoted by $\Delta_1$, $\Delta_2$, $\Delta_3$, $\Delta_4$ respectively. Suffices are used for key points in these triangles. Thus $H_k$ is the orthocentre of $\Delta_k$, $N_k$ is the nine-point centre of $\Delta_k$ and $G_k$ is the centroid of $\Delta_k$, $k = 1, 2, 3, 4$. The vector positions corresponding to the points are $\mathbf{h_1} = \mathbf{b} + \mathbf{c} + \mathbf{d}$, $\mathbf{h_2} = \mathbf{a} + \mathbf{c} + \mathbf{d}$, $\mathbf{h_3} = \mathbf{a} + \mathbf{b} + \mathbf{d}$, $\mathbf{h_4} = \mathbf{a} + \mathbf{b} + \mathbf{c}$ and, as in any triangle, $\mathbf{n_k} = \frac{1}{2}\mathbf{h_k}$ and $\mathbf{g_k} = \frac{1}{3}\mathbf{h_k}$. The following facts are now easy to establish:

- (i) $AH_1$, $BH_2$, $CH_3$, $DH_4$ are concurrent at the point *N* with vector position $\frac{1}{2}(\mathbf{a} + \mathbf{b} + \mathbf{c} + \mathbf{d})$;
- (ii) $AN_1$, $BN_2$, $CN_3$, $DN_4$ are concurrent at the point *G* with vector position $\frac{1}{3}(\mathbf{a} + \mathbf{b} + \mathbf{c} + \mathbf{d})$;
- (iii) $AG_1$, $BG_2$, $CG_3$, $DG_4$ are concurrent at a point *M* with vector position $\frac{1}{4}(\mathbf{a} + \mathbf{b} + \mathbf{c} + \mathbf{d})$;
- (iv) $H_1$, $H_2$, $H_3$, $H_4$ lie on a circle of radius *R* with centre *H*, which has position vector $(\mathbf{a} + \mathbf{b} + \mathbf{c} + \mathbf{d})$ and moreover $H_1H_2H_3H_4$ is the image of *ABCD* under a 180° rotation about *N*;
- (v) $N_1$, $N_2$, $N_3$, $N_4$ lie on a circle of radius $\frac{1}{2}R$ with centre *N* and $N_1N_2N_3N_4$ is homothetic with $H_1H_2H_3H_4$, with centre *O* and enlargement factor ½;
- (vi) The nine-point circles of triangles $\Delta_1$, $\Delta_2$, $\Delta_3$, $\Delta_4$ , each having radius $\frac{1}{2}R$ therefore all pass through *N*;
- (vii) A known result is that a rectangular hyperbola through the vertices of a triangle passes through its orthocentre. It follows that the rectangular hyperbola $\Sigma$ through *A*, *B*, *C*, *D* passes through all of $H_1$, $H_2$, $H_3$, $H_4$;
- (viii) Another known result is that the centre of a rectangular hyperbola passing through the vertices of a triangle lies on its nine-point circle. In consequence of (vii) the centre of $\Sigma$ is the point *N*;
- (ix) $G_1$, $G_2$, $G_3$, $G_4$ lie on a circle of radius $\frac{1}{3}R$ with centre *G*, which has position vector $\frac{1}{3}(\mathbf{a} + \mathbf{b} + \mathbf{c} + \mathbf{d})$;
- (x) The points *O*, *M*, *G*, *N*, *H* are collinear and if *O* and *H* are given co-ordinates 0, 1 on this line, then *M*, *G*, *N* have co-ordinates ¼, ⅓, ½ respectively.

See Fig.1 for an illustration of these properties.

## 2. Setting the scene

It is clear from Speckman's work on indirect similar perspective triangles that if we now choose any point $P$ on the rectangular hyperbola defined in Section 1 Result (vii), it is possible to find the Hagge circle of $P$ with respect to each of the four triangles $\Delta_1$, $\Delta_2$, $\Delta_3$, $\Delta_4$ and to draw them on the same diagram. These are the four Hagge circles $\Gamma_1$, $\Gamma_2$, $\Gamma_3$, $\Gamma_4$ of the article heading and they are shown in Fig. 2. Note that as the axes of inverse similarity through $P$ are parallel to the asymptotes of the rectangular hyperbola $\Sigma$, these axes coincide for each of the Hagge circles.

A consistent notation is essential and this we now describe. See Fig. 2. The centres of the circles are denoted by $Q_k$, $k = 1, 2, 3, 4$. For reasons that become clear shortly we re-label the four orthocentres $A_1$, $B_2$, $C_3$, $D_4$, so that, for example $D_4$ is the orthocentre of triangle $ABC$. The lines $AP$, $BP$, $CP$, $DP$ meet the circumcircle $\Gamma$ at points $A'$, $B'$, $C'$, $D'$. We now consider the labelling of points on the Hagge circle $\Gamma_4$, which is the Hagge circle of $P$ with respect to triangle $ABC$. They all carry the subscript 4. The reflection of $A'$ in $BC$ we denote by $A_4'$, the reflection of $B'$ in $CA$ we denote by $B_4'$ and the reflection of $C'$ in $AB$ we denote by $C_4'$. These replace the labels $U$, $V$, $W$ used in Article: CJB/2010/1 [1]. The points where $A_4'P$, $B_4'P$, $C_4'P$, $D_4P$ meet $\Gamma_4$ are denoted by $A_4$, $B_4$, $C_4$, $D_4'$ respectively. The first three of these replace the labels $X$, $Y$, $Z$ in Article 1. Points with suffices 1, 2, 3 are similarly defined on circles $\Gamma_1$, $\Gamma_2$, $\Gamma_3$. The labels $Pg_1$, $Pg_2$, $Pg_3$, $Pg_4$ are given to the isogonal conjugates of $P$ with respect to triangles $\Delta_1$, $\Delta_2$, $\Delta_3$, $\Delta_4$ respectively. The double lines of inverse symmetry through $P$ are denoted by $L$ and $L'$. The centre of the rectangular hyperbola $\Sigma$ is denoted by $M$.

In Section 3 we use Cartesian co-ordinates, origin $M$ and the asymptote parallel to $L$ is chosen as the $x$-axis. In this way $\Sigma$ has equation $xy = 1$ (by choice of scale), and we use a parameter $t$ on the hyperbola so that its points have co-ordinates $(t, 1/t)$. The points $A$, $B$, $C$, $D$, $P$ are given parameters $a$, $b$, $c$, $d$, $p$, where it is known, since $ABCD$ is cyclic, that $abcd = 1$ and the parameters of $A_1$, $B_2$, $C_3$, $D_4$, since they are the orthocentres, have parameters $-a$, $-b$, $-c$, $-d$ respectively.

The following theorems now hold:

**Theorem 1**
The centres, $Q_1$, $Q_2$, $Q_3$, $Q_4$, of the four Hagge circles are collinear.

**Theorem 2**
The points $A_1$, $A_2$, $A_3$, $A_4$, $A_1'$, $A_2'$, $A_3'$, $A_4'$ are collinear and similarly for the points $B_k$, $B_k'$, $k = 1, 2, 3, 4$ and $C_k$, $C_k'$, $k = 1, 2, 3, 4$ and $D_k$, $D_k'$, $k = 1, 2, 3, 4$.

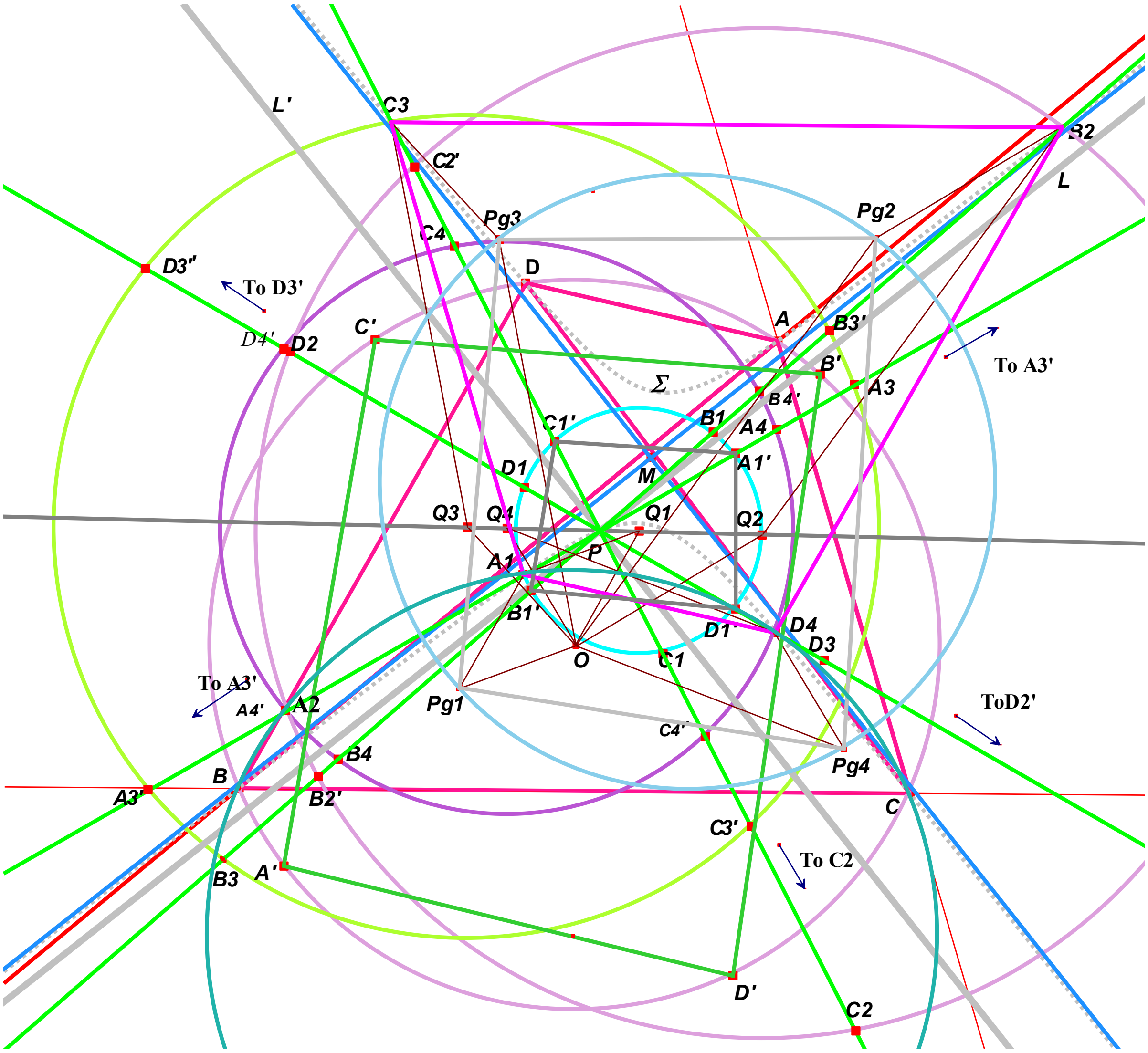

**Fig.2**
**Illustrating the notation and results for the four Hagge circles**

**Theorem 3**
The quadrilateral $Pg_1Pg_2Pg_3Pg_4$ is similar to the quadrilateral $A'B'C'D'$.

## 3. The analysis

**The co-ordinates of *O***

The chord *AB* has equation

$$aby + x = a + b. \qquad (3.1)$$

The midpoint of *AB* has co-ordinates (½(*a* + *b*), ½(*a* + *b*)/(*ab*)). It follows that the perpendicular bisector of *AB* has equati0n

$$2ab(y - abx) = (a + b)(1 - a^2b^2). \qquad (3.2)$$

The perpendicular bisector of *BC* has equation

$$2bc(y - bcx) = (b + c)(1 - b^2c^2). \qquad (3.3)$$

These meet at *O*, the centre of *Γ*, at the point with co-ordinates (½(*a* + *b* + *c* + *d*), ½(1/*a* + 1/*b* + 1/*c* + 1/*d*). Here we have used *abcd* = 1.

It can be shown that the radius of the circumcircle is

$$\tfrac{1}{2}\sqrt{(a^2 + b^2 + c^2 + d^2 + 1/a^2 + 1/b^2 + 1/c^2 + 1/d^2)}. \qquad (3.4)$$

**The equation of the circumcircle *Γ***

From the above information this is easily proved to be

$$2x^2 - 2x(a + b + c + d) + 2y^2 - 2y(1/a + 1/b + 1/c + 1/d) + (ab + ac + ad + bc + bd + cd + 1/ab + 1/ac + 1/ad + 1/bc + 1/bd + 1/cd) = 0. \qquad (3.5)$$

Alternatively, in terms of *a*, *b*, *c* only, it has the form

$$abcx^2 - (abc(a + b + c) + 1)x + abcy^2 - (a^2b^2c^2 + ab + bc + ca)y + (abc)(ab + bc + ca) + a + b + c = 0. \qquad (3.6)$$

**A digression**

The reflection of the point with co-ordinates (*f*, g) in the line with equation *lx* + *my* =*n* is the point with co-ordinates (*h*, *k*) where

$$h = \{f(m^2 - l^2) + 2l(n - gm)\}/(l$$
$$k = \{2mn + g(l^2 - m^2) - 2flm\}/(l^2 + m^2).$$

This is straightforward and is left to the reader. This is used later to reflect $ABCD$ in the line $L$ and then dilate through $P$ to obtain the Hagge circle $A_4B_4C_4D_4$.

**The indirect similarity**

We now determine the equations governing the indirect similarity which maps the circumcircle into the Hagge circle of $P$ with respect to $\Delta_4$. This is done by following what happens to the point $D(d, 1/d)$. We know the image of this point is the orthocentre $D_4(-1/abc, -abc)$ and we know the mapping is effected by a reflection in the line $L$, with equation $y = 1/p$ followed by an enlargement (reduction) by a factor $PD_4/PD$. As far as the reflection is concerned we may use the analysis of the last paragraph with $l = 0$, $m = 1$, $n = 1/p$. The result of this reflection on $D$ is to produce the point with co-ordinates $(d, 2/p - 1/d)$.

The enlargement (reduction) factor (comparing $x$-co-ordinates of the points concerned) is equal to $(p + d)/(p - d) = (abcp + 1)/(abcp - 1)$. It may now be checked that the reflection in $L$ followed by the dilation with this enlargement (reduction) factor takes the point with co-ordinates $(h, k)$ to the point with co-ordinates $(x, y)$ where

$$x = (abchp - 2p + h)/(abcp - 1), \qquad (3.7)$$

and

$$y = (abc(kp - 2) + k)/(1 - abcp). \qquad (3.8)$$

Using $h = a$, $k = 1/a$ we deduce the co-ordinates of $A_4$ to be

$$((a^2bcp + a - 2p)/(abcp - 1),\ (2a^2bc - abcp - 1)/(a(abcp - 1))),$$

with similar expressions for the co-ordinates of $B_4$ and $C_4$ by using $b$ and $c$, in $x$ and $y$ above, instead of $a$.

**The equations of the Hagge circles**

From here it is quite an elaborate calculation to obtain the equation of the Hagge circle $A_4B_4C_4$ and check the fundamental theorem that $D_4$ lies on this circle. The computer algebra package *DERIVE* was used to perform the calculations, and the result is that $\Gamma_4$ has equation

$$\begin{aligned} &abc(abcp - 1)(x^2 + y^2) - (a^2b^2c^2p(a + b + c) + abc(a + b + c) - 3abcp + 1)x \\ &\quad + (a^3b^3c^3p - 3a^2b^2c^2 + (abcp + 1)(bc + ca + ab))y \\ &- 2a^3b^3c^3 + a^2b^2c^2p(a + b + c) + abc(bc + ca + ab) - (abcp + 1)(a + b + c) - 2p = 0. \end{aligned} \qquad (3.9)$$

The equations of the other Hagge circles follow immediately by using other triplets of parameters instead of $a$, $b$, $c$. The co-ordinates of the centres of these circles may now be written down and it has been checked, using *DERIVE*, that any three of the four centres are collinear. Theorem 1 is now proved, but it also follows from a geometrical argument. The centre $O$ of the circumcircle is the common circumcentre of all four triangles $\Delta_1$, $\Delta_2$, $\Delta_3$, $\Delta_4$. If we now reflect $O$ in the two lines of inverse similarity, it maps into a pair of points collinear with $P$. Since the centres of the four Hagge circles are now the images of these points by means of dilations through $P$ with enlargements (reductions) using different factors, the resulting images all lie on this line, and this line also passes through $P$. Theorem 2 follows by a similar argument, bearing in mind that we already know primed and unprimed pairs of points such as $A_4$, $A_4'$ are collinear with $P$.

**The quadrilateral $A'B'C'D'$**

The equation of the line $AP$ is

$$x + apy = a + p. \qquad (3.10)$$

This meets the circumcircle again at the point $A'$ with co-ordinates $(x, y)$ where

$$x = \{ap^2(abc(b + c) + 1) - p(a^2b^2c^2 + a(b + c) - bc) + abc\}/\{bc(a^2p^2 + 1)\}, \qquad (3.11)$$

$$y = \{abcp^2 + p(a^2bc - abc(b + c) - 1) + ab^2c^2 + b + c\}/\{bc(a^2p^2 + 1)\}. \qquad (3.12)$$

The point $D'$ has co-ordinates similar to these, but with $d$ replacing $a$.
From these co-ordinates we can work out $(A'D')^2$ and the result is

$$\{(a - d)^2(b - p)^2(c - p)^2(b^2c^2 + 1)\}/\{b^2c^2(a^2p^2 + 1)(d^2p^2 + 1)\}. \qquad (3.13)$$

**The quadrilateral $Pg_1Pg_2Pg_3Pg_4$**

The co-ordinates $(h, k)$ of $Pg_4$, the isogonal conjugate of $P$ with respect to triangle $ABC$ are given by

$$h = (a + b + c - p)/(1 - abcp), \qquad (3.14)$$

$$k = (p(bc + ca + ab) - abc)/(abcp - 1). \qquad (3.15)$$

This may be checked as follows: Let the line from $A$ to $(h, k)$ meet $\Gamma$ at $A''$, then it is easy to show that $A'A''$ is parallel to $BC$. The symmetry of $h$, $k$ with respect to $a$, $b$, $c$ now proves that $(h, k)$ are the co-ordinates of $Pg_4$, as this calculation reflects one of the standard constructions for an isogonal conjugate. Indeed the parallel $A'A''$ to $BC$ indicates that $AP$ and $APg_4$ are reflections of each other in the internal bisector of angle $A$. The co-ordinates of $Pg_1$ now follow from $(h, k)$ by exchanging $a$ and $d$.
From these co-ordinates we can work out $(Pg_1Pg_4)^2$ and the result is

$$\{(a - d)^2(b^2c^2 + 1)(b^2p^2 + 1)(c^2p^2 + 1)\}/\{(abcp - 1)^2(bcdp - 1)^2\}. \qquad (3.16)$$

**The similarity of the quadrangles**

The ratio of the squares of the side lengths and diagonals of the quadrangles $Pg_1Pg_2Pg_3Pg_4$ and $A'B'C'D'$ may now be calculated and is the totally symmetric expression

$$\{(a^2p^2 + 1)(b^2p^2 + 1)(c^2p^2 + 1)(d^2p^2 + 1)\}/\{(a - p)^2(b - p)^2(c - p)^2(d - p)^2\}. \qquad (3.17)$$

Flat 4,
Terrill Court,
12-14 Apsley Road,
BRISTOL BS8 2SP